\documentclass[a4paper]{article}

\usepackage{amsfonts}
\usepackage{amsmath}
\usepackage{amssymb}
\usepackage{graphicx}
\usepackage{mathtools}

\DeclarePairedDelimiter\abs{\lvert}{\rvert}%
\DeclarePairedDelimiter\norm{\lVert}{\rVert}%

\makeatletter
\let\oldabs\abs
\def\abs{\@ifstar{\oldabs}{\oldabs*}}
\let\oldnorm\norm
\def\norm{\@ifstar{\oldnorm}{\oldnorm*}}
\makeatother

\def\||{{\, || \,}}

\def\bor{\mathop{\mathord{\lor}\!\!\!\raise4pt\hbox{$\scriptscriptstyle 2$}\,}}
\def\band{\mathop{\mathord{\land}\!\!\!\lower2pt\hbox{$\scriptscriptstyle 2$}\,}}
\def\OR{{\tt OR}}

\newtheorem{theorem}{Theorem}

\def\mydot{ \mathrel{\raisebox{2pt}{\scalebox{0.5}{\textbullet}}}}

\begin{document}
\title{Why Square Roots of Probabilities?}
\author{Kevin H. Knuth\\ Departments of Physics and Informatics\\
University at Albany (SUNY)\\ Albany NY 12222, USA}

\date{\today}
\maketitle

\begin{abstract}
Square roots of probabilities appear in several contexts, which suggests that they are somehow more fundamental than probabilities.  Square roots of probabilities appear in expressions of the Fisher-Rao Metric and the Hellinger-Bhattacharyya distance.  They also come into play in Quantum Mechanics via the Born rule where probabilities are found by taking the squared modulus of the quantum amplitude.  Why should this be the case  and why do these square roots not arise in the various formulations of probability theory?

In this short, inconclusive exploration, I consider quantifying a logical statement with a vector defined by a set of components each quantifying one of the atomic statements defining the hypothesis space.  I show that conditional probabilities (bi-valuations), such as $P(x|y)$, can be written as the dot product of the two vectors quantifying the logical statements $x$ and $y$ each normalized with respect to the vector quantifying the conditional $y$.  The components of the vectors are proportional to the square root of the probability.  As a result, this formulation is shown to be consistent with a concept of orthogonality applied to the set of mutually exclusive atomic statements such that the sum rule is represented as the sum of the squares of the square roots of probability.
\end{abstract}

\medskip\noindent{Key Words: }
{probability, orthogonality, information, information geometry, quantum mechanics, Fisher-Rao metric, Hellinger-Bhattacharyya distance}

\section{Introduction}
Square roots of probabilities appear in a variety of interesting contexts ranging from information geometry and statistical manifolds in the form of the Fisher-Rao metric
\begin{equation}
g_{ij}(\theta) = 4 \int{ \frac{\partial \sqrt{p(x|\theta)}}{\partial \theta_i} \frac{\partial \sqrt{p(x|\theta)}}{\partial \theta_j} \, dx }
\end{equation}
and the Hellinger-Bhattacharyya Distance
\begin{equation}
H^2(P,Q) = \frac{1}{2} \int{ \left( \sqrt{\frac{dP}{d\lambda}} - \sqrt{\frac{dQ}{d\lambda}} \right)^2 d\lambda }
\end{equation}
to quantum mechanics via the Born rule where probabilities are found by taking the squared modulus of complex-valued quantum amplitudes
\begin{equation}
p(x) = | \psi(x) |^2.
\end{equation}
While quantum amplitudes are clearly fundamental, and in some sense foundationally understood in terms of symmetries \cite{GKS:PRA}\cite{GK:Symmetry}, the roles that square roots (or roots) of probabilities play in various applications \cite{Daum:1986}\cite{Brigo+etal:1998}\cite{Chizat+etal:2015}\cite{Pfaff+etal:2015} have an unclear and less well-understood\footnote{At the very least, this is certainly less well-understood by me.} foundational significance.  This includes the area where quantum mechanics and information geometry intersect \cite{Goyal:2008}.

This paper represents a short, inconclusive exploration into the fact that square roots of probabilities appear to be somehow fundamental, but do not appear to play a role in any of the major foundational approaches to probability theory \cite{DeFinetti:1931}\cite{Kolmogorov:1950}\cite{Cox:1946}\cite{Knuth&Skilling:2012}.

\section{Quantification with Vectors}
Let me propose---out of the blue---that we quantify statements with a vector.  This is not a derivation---it is instead an exploration.

Consider a hypothesis space based on three atomic (mutually exclusive) statements $a_1$, $a_2$, and $a_3$, and let us quantify each statement with a vector, such that
the statement $a_1$ is quantified by the vector $(q(a_1), 0, 0)$, where $q$ is a function that takes its argument, which is a statement, to a real number.  We then have the following:
\begin{align}
a_1 \quad &\mbox{is quantified by} \quad v(a_1) := (q(a_1), 0, 0)  \nonumber \\
a_2 \quad &\mbox{is quantified by} \quad v(a_2) := (0, q(a_2), 0)  \nonumber \\
a_3 \quad &\mbox{is quantified by} \quad v(a_3) := (0, 0, q(a_3)).  \nonumber
\end{align}

The join, or disjunction, of two statements, denoted by $\vee$, is a statement found by applying the logical \OR.  As such, the join operation is associative and commutative. Since we want to ultimately rank the statements, given the properties of the join operation we can expect
\footnote{Again, this is not a derivation.  However, Skilling has an unpublished derivation that could be employed to make this assertion of additivity rigorous.}
that one should sum the vectors of atomic statements: \cite{Knuth:FQXI2015}
\begin{equation}
a_1 \vee a_2 \quad \mbox{is quantified by} \quad v(a_1 \vee a_2) := v(a_1) + v(a_2) = (q(a_1), q(a_2), 0).  \nonumber \\
\end{equation}
Vectors representing non-atomic statements should be summed using the inclusion-exclusion formula to avoid overcounting \cite{Knuth:FQXI2015}.

\section{Bi-Valuations and Inner Products}
Let us consider bi-valuations, such as $P(x|y)$, to represent a projection from one statement $x$ onto another statement $y$, which is called the context, defined by a context-normalized inner product where
\begin{equation}
P(x | y) := v(x) \mydot \frac{v(y)}{|v(y)|^2}.
\end{equation}
Note that the fact that the inner product is normalized with respect to the squared magnitude of the vector representing the context (second argument) introduces an important asymmetry.

Next we show, via several examples, that this forces the vector components to be square roots of probabilities.
For example, since the atoms $a_1$, $a_2$, and $a_3$ are assumed to be mutually exclusive and exhaustive, the join of all three atoms is the truism.  This means that the probability of any one of the atoms can be expressed as the conditional probability (bi-valuation) $p_i = P(a_i | a_1 \vee a_2 \vee a_3)$. For $a_i = a_1$, we then have that:
\begin{align}
P(a_1 | a_1 \vee a_2 \vee a_3) &= v(a_1) \mydot \frac{v(a_1 \vee a_2 \vee a_3)}{|v(a_1 \vee a_2 \vee a_3)|^2}\\
&= (q(a_1), 0, 0) \mydot \frac{(q(a_1), q(a_2), q(a_3))}{q(a_1)^2 + q(a_2)^2 + q(a_3)^2}\\
&= \frac{q(a_1)^2}{q(a_1)^2 + q(a_2)^2 + q(a_3)^2}\\
&= \frac{m(a_1)}{m(a_1)+m(a_2)+m(a_3)},
\end{align}
where $m(a_1) + m(a_2) + m(a_3) = 1$ and $m(a_i) = q(a_i)^2$ is proportional to the probability of $a_1$.  This implies that the vector components $q(a_i)$ are square roots of probabilities: $q(a_i) = \sqrt{p_i}$.

To ensure that this concept of inner products is consistent with probability theory, we can consider switching the arguments of the bi-valuation above.  This results in
\begin{align}
P(a_1 \vee a_2 \vee a_3 | a_1) &= v(a_1 \vee a_2 \vee a_3) \mydot \frac{v(a_1)}{|v(a_1)|^2}\\
&= (q(a_1), q(a_2), q(a_3)) \mydot \frac{(q(a_1),0,0)}{q(a_1)^2}\\
&= \frac{q(a_1)^2}{q(a_1)^2}\\
&= 1,
\end{align}
which is as expected since $a_1$ implies $a_1 \vee a_2 \vee a_3$.

Let us now consider two statements, $a_1$ and $a_2$, which are mutually exclusive.  Note that their vectors are orthogonal, which results in a zero projection
\begin{align}
P(a_1 | a_2) &= v(a_1) \mydot \frac{v(a_2)}{|v(a_2)|^2}\\
&= (q(a_1), 0,0) \mydot \frac{(0,q(a_2),0)}{q(a_2)^2}\\
&= \frac{0}{q(a_2)^2}\\
&= 0.
\end{align}

Last, let us consider a more limited context
\begin{align}
P(a_1 | a_1 \vee a_2) &= v(a_1) \mydot \frac{v(a_1 \vee a_2)}{|v(a_1 \vee a_2)|^2}\\
&= (q(a_1), 0, 0) \mydot \frac{(q(a_1), q(a_2),0)}{q(a_1)^2 + q(a_2)^2}\\
&= \frac{q(a_1)^2}{q(a_1)^2 + q(a_2)^2}\\
&= \frac{m(a_1)}{m(a_1)+m(a_2)},
\end{align}
which is as expected.


\section{Determinants and Volume}
Packing the vectors quantifying the atomic statements into a matrix
$$
g =
\begin{pmatrix}
q(a_1) & 0 & 0 \\
0 & q(a_2) & 0 \\
0 & 0 & q(a_3)
\end{pmatrix}
$$
and taking the determinant
$$
g =
\begin{vmatrix}
q(a_1) & 0 & 0 \\
0 & q(a_2) & 0 \\
0 & 0 & q(a_3)
\end{vmatrix}
$$
yields a volume given by
\begin{align}
det(g) &= q(a_1)q(a_2)q(a_3) \nonumber \\
&= \sqrt{p_1 p_2 p_3}.
\end{align}
Assigning a prior probability based on the inverse volume\footnote{This is true regardless of whether such an assignment would be appropriate.} yields the non-informative prior for the multinomial. \cite{Box&Tiao:1992}

\section{Conclusions}
Statements are quantified by vectors with a dimensionality given by the number of atomic statements.  The vectors quantifying the atomic statements form a basis, such that each of the atomic vectors has one and only one non-zero component proportional to the square root of the probability of that statement.  In general, mutually exclusive statements are quantified by vectors that are orthogonal to one another.  The probability of one statement given another is proportional to the inner product of the two vectors.  This results in a real-valued Hilbert space.

In general, entities that are orthogonal to one another are often quantified by quantities whose squares sum when the entities are joined.  Azc\'el and Dhombres \cite{aczel+dhombres:1989} present the following theorem which establishes the precise relation
\begin{theorem}[Azc\'el and Dhombres: Theorem 11.17]
Let $(E,<,>)$ be a real inner product space of dimension at least 2. A continuous function $f:E \rightarrow \mathbb{R}$ is orthogonally additive if, and only if, there exists a real constant $a$ and a continuous linear function $h:E \rightarrow \mathbb{R}$ such that
\begin{equation}
f(x) = a \norm{x}^2 + h(x) \; (x \in E).
\end{equation}
\end{theorem}
While this exposition does not constitute a rigorous treatment, it appears that since we can represent statements as vectors with a well-defined inner product, this results in conditional probabilities, which quantify these vectors in a way (represented by the function $f$), that depends, at least in part on the square modulus of the vector.  Here the sum rule of probability in the case of two mutually exclusive statements becomes a sum of squares of square roots of probabilities as one would expect for orthogonal entities:
\begin{align}
P(a_1 \vee a_2 | a_1 \vee a_2 \vee a_3) &= P(a_1 | a_1 \vee a_2 \vee a_3) + P(a_2 | a_1 \vee a_2 \vee a_3)\\
&= \frac{q(a_1)^2}{Z} + \frac{q(a_2)^2}{Z}\\
&= \frac{q(a_1)^2 + q(a_2)^2}{Z},
\end{align}
where $Z = q(a_1)^2 + q(a_2)^2 + q(a_3)^2$.
So it could be that the square roots of probabilities are fundamental in the sense that they allow one to quantify the combination of orthogonal entities via the sum of squares.

Such a formulation of probability theory may promise to resolve some of these mysteries surrounding square roots of probability, tie in more closely with some of the geometric concepts discussed here and elsewhere, and perhaps shed more light on quantum mechanics and the fact that the Born rule involves taking the square of the modulus of the quantum amplitude.

\section{Acknowledgements}
I would like to thank Ariel Caticha and James Walsh for their careful reading of this manuscript.  I would also like to thank Fredrick Daum for bringing several of the applications of square roots of probabilities and densities to my attention.

\bibliographystyle{amsplain}
\bibliography{C:/Users/KK952431/kevin/files/papers/bibliography/knuth}

\end{document}